\documentclass{article}
\usepackage{amssymb}

\newtheorem{theorem}{Theorem}[section]
\newtheorem{proposition}[theorem]{Proposition}
\newtheorem{conjecture}[theorem]{Conjecture}
\newtheorem{procedure}[theorem]{Procedure}
\newtheorem{lemma}[theorem]{Lemma}
\newtheorem{example}[theorem]{Example}
\newtheorem{definition}[theorem]{Definition}
\newtheorem{remark}[theorem]{Remark}
\newtheorem{slogan}[theorem]{Slogan}

\newcommand{\dP}{{\mathbb P}}
\newcommand{\dC}{{\mathbb C}}
\newcommand{\dZ}{{\mathbb Z}}
\newcommand{\dQ}{{\mathbb Q}}
\newcommand{\dR}{{\mathbb R}}

\newcommand{\X}{{\mathcal X}}
\newcommand{\Y}{{\mathcal Y}}
\newcommand{\D}{{\mathcal D}}
\newcommand{\F}{{\mathcal F}}
\newcommand{\M}{{\mathcal M}}
\newcommand{\sheafO}{{\mathcal O}}
\newcommand{\aut}{\mathop{\rm Aut\,}\nolimits}
\newcommand{\pic}{\mathop{\rm Pic\,}\nolimits}
\newcommand{\diff}{\mathop{\rm Diff}\nolimits}
\newcommand{\mir}{\mathop{\rm Mir}\nolimits}
\newcommand{\db}[1]{D^b(#1)}
\newcommand{\ethrm}{\hspace*{\fill} $\square$ \vspace{.1in}}
\def\mapdownright#1{\Big\downarrow\rlap{$\vcenter{\hbox{$\scriptstyle#1$}}$}}
\def\mapdownleft#1{\rlap{$\vcenter{\hbox{$\scriptstyle#1$}}$}\,\,\,\Big\downarrow} 

\begin{document}
 
\begin{center}
{\Large Diffeomorphisms and families of Fourier--Mukai

transforms in mirror symmetry}

\vspace{0.2in}

{\large Bal\' azs Szendr\H oi}

\vspace{0.15in}

{\it Mathematics Institute, University of Warwick {\rm and} Alfr\'ed R\'enyi 

Institute of Mathematics, Hungarian Academy of Sciences}

\vspace{0.2in}

\end{center}

\begin{abstract}
Assuming the standard framework of mirror symmetry, a conjecture is formulated describing how the diffeomorphism group of a Calabi--Yau manifold $Y$ should act by families of Fourier--Mukai transforms over the complex moduli space of the mirror $X$. The conjecture generalizes a proposal of Kontsevich relating monodromy transformations and self-equivalences. Supporting evidence is given in the case of elliptic curves, lattice-polarized K3 surfaces and Calabi--Yau threefolds. A relation to the global Torelli problem is discussed.
\end{abstract}

\section*{Introduction} 

Derived categories of coherent sheaves entered the mirror symmetry scene 
with the paper~\cite{kontsevich}. The mirror relationship as envisaged by 
Kontsevich is an equivalence of $A_\infty$-categories built out of two 
manifolds $X$ and $Y$, both equipped with (complexified) 
symplectic forms and complex Calabi--Yau structures.
One of these categories is $\db{X_t}$, the bounded derived category 
of coherent sheaves on the complex manifold $X_t$, 
or a twisted version thereof. The other category is the derived Fukaya category
${D^b{\rm Fuk}}(Y,\omega_0)$ of the symplectic manifold $(Y, \omega_0)$, 
a category constructed from Lagrangian submanifolds 
and local systems on them, with morphism spaces given by Floer homology. 

The origin of the equivalence of the two categories was not specified by 
Kontsevich. The picture was filled in by a proposal of Strominger--
Yau--Zaslow~\cite{syz}, based on 
arguments coming from non-per\-tur\-ba\-tive string theory.
According to~\cite{syz}, mirror manifolds should be fibered into
middle-dimensional real special Lagrangian
tori over a common base space; moreover, 
the two fibrations should in a suitable sense be
dual to each other. As it was subsequently realized, 
this should supply the equivalence of categories proposed by Kontsevich 
by an analytic form of the Fourier--Mukai transform, 
converting Lagrangian submanifolds equipped with a local system on one 
manifold into holomorphic bundles or more generally sheaves on the 
mirror; see e.g.~\cite{slagtohym}.

The next question that arises is the origin of the torus fibrations; 
recently, some proposals have been put forward by 
Gross--Wilson~\cite{gross_wilson2} and 
Kontsevich--Soibelman~\cite{konts_soi} 
to the effect that these should arise from certain 
degenerations of the complex structure on the Calabi--Yau manifolds. These
degenerations, to so-called maximally degenerate boundary points or 
complex cusps in the Calabi--Yau moduli space, have been known to play 
a major role in mirror symmetry since the heroic age~\cite{cogp}; 
the ideas of~\cite{gross_wilson2} and \cite{konts_soi} 
show their significance in the new order of things.

Many of the details of the torus fibrations, the definitions of the 
categories and their equivalence, and the relation to more traditional 
notions such as curve-counting
generation functions and variations of Hodge structures are missing or at best
conjectural. However, this should not necessarily prevent one from 
investigating
some further issues such as symmetries on the two sides of the mirror 
map. The study of these issues was also initiated by Kontsevich. 
He pointed out that the categorical equivalence implies a correspondence 
between symplectomorphisms on a Calabi--Yau manifold and self-equivalences of 
the derived category of coherent sheaves of its mirror, and gave some explicit 
examples of this relationship. In particular, he related symplectomorphisms 
arising from certain 
monodromy transformations to certain relatively simple self-equivalences 
on the mirror side in the case of the quintic, 
and proposed some more general constructions of 
self-equivalences. These ideas were carried further by Seidel and 
Thomas~\cite{st}
and in a different direction by Horja~\cite{horja},~\cite{horja2}. 
Aspinwall~\cite{aspinwall} and 
Bj\"orn, Curio, Hern\'andez Ruip\'erez and Yau~\cite{bchy} studied
the relation between monodromies and self-equivalences
in certain concrete cases. 

A common feature of these works is that they present a static picture: 
the symplectic form and complex structure are frozen, and symmetries compared. 
The main contribution of the present work is that by investigating 
these notions in families, some additional insight can be gained.
 
The language of kernels (derived correspondences) is used to define
a group that I term the categorical mapping group of a family of
complex manifolds. The main conjecture of the paper states that 
for a pair $(X,Y)$ of mirror Calabi--Yau manifolds, there should be a  
homomorphism from the diffeomorphism group of $Y$ to the categorical 
mapping group of the total space of the complex moduli space of $X$. 
A weaker version of the conjecture concerns cohomology actions on the two 
sides of the mirror map. 

I investigate the proposed conjecture for families of Calabi--Yau 
manifolds of dimension at most three. After a brief look at the elliptic
curve case, where the relevant group action was already known to 
Mukai\footnote{The reader will excuse me for a short historical digression. 
Mukai remarks in the introduction of~\cite{mukai_abelian} that 
the result that the derived category of an elliptic curve carries an action of 
$SL(2,\dZ)$ modulo the shift `seems to be significant'. 
This must be one of the first explicit hints to mirror symmetry in the 
mathematics literature. Another early hint is 
discussed in~\cite[Section 8]{morrison_compactifications}: 
a picture of Mori about the cone of curves of an abelian surface closely 
resembles a picture of Mumford about the compactification of a Hilbert modular 
surface; these cones become duals 
under mirror symmetry. It is a curious fact that Mori drew his picture 
during the fall of 1979 which is exactly the time 
when~\cite{mukai_abelian} was submitted.}, 
I consider the case of K3 surfaces.  Derived categories of sheaves on K3s
have only been studied in the projective case, 
to which I restrict; this can be done elegantly in families 
by considering lattice-polarized K3 surfaces, a notion due to Nikulin.
The diffeomorphism group has to be restricted also; 
in effect, the lattice polarization partitions 
$H^2$ (generically) into algebraic and transcendental 
parts a priori, and the allowable diffeomorphisms are supposed to respect this
choice. One further issue that arises is an analogue of a theorem of Donaldson
on the complex side, regarding orientation on the cohomology of K3 surfaces.
Unfortunately, I cannot prove the required statement at present; it is 
formulated as Conjecture~\ref{orient}. Assuming this, 
Theorem~\ref{K3_main_theorem} confirms that for any allowable diffeomorphism, 
there is a family of categorical equivalences 
with the correct cohomology action. As an illustration I investigate this 
correspondence in some detail for toric families. 

In the threefold case, I only discuss examples of elements of the
categorical mapping group, mostly coming from birational contractions
on the threefolds. A typical case of Horja's work emerges in a novel 
fashion, showing that indeed, the present version of the conjecture is needed
to get a full picture. To conclude, I point out a relation to 
(counterexamples to) the global Torelli problem. 

The purpose of this paper is to discuss the proposed framework. Detailed
proofs of the assertions, together with some further examples and 
applications, will be given elsewhere. 

\subsection*{Acknowledgments}
Conversations and correspondence with Tom Brid\-ge\-land, 
Andrei C\u al\-d\u a\-ra\-ru, Mark Gross, Paul Seidel and Richard Thomas 
during the course of this work were very helpful indeed.  
I thank the MPI for hospitality during a brief stay in September 2000, 
where some of these ideas took shape. I especially thank Yuri Manin for the 
warm welcome and stimulating conversations; he also provided the starting 
point for these ideas with questions in~\cite{manin}. 

\section{Basic definitions} 

For the purposes of the paper, a {\it Calabi--Yau manifold} is 
complex manifold 
$X_0$ with holonomy exactly $SU(n)$, in particular trivial canonical bundle. 
I assume everywhere below that $n=\dim_\dC(X_0)\leq 3$. The subscript means
that a particular complex structure is chosen on the differentiable 
manifold $X$; sometimes I write $X_0$ as a pair $(X, I_0)$ to make the 
complex structure explicit in notation. $\db{X_0}$ or $\db{X, I_0}$
denote the 
bounded derived category of coherent sheaves on $X_0$. The {\it Mukai map} 
\[ v\colon \db{X_0}\rightarrow H^*(X_0, \dC)
\]
is defined for $U\in\db{X_0}$ by
\[
v(U)={\rm ch}(U)\cup\sqrt{{\rm td}(X_0)}. 
\]

The {\it Mukai pairing} on cohomology 
\[ 
H^*(X,\dZ) \times  H^*(X,\dZ)  \rightarrow \dZ+i\dZ  
\]
is defined by 
\[(\alpha_0+ \ldots +\alpha_{2n}) \cdot (\beta_0+\ldots +\beta_{2n}) = (-1)^{n-1}\sum_{j=0}^{2n} i^j \int_Y \alpha_j \cup \beta_{2n-j}
\]
where $\alpha_i,\beta_i\in H^i(X, \dZ)$; this pairing extends linearly to 
rational and complex cohomology. Let $\tilde L=H^*(X,\dZ)$ be 
the $\dZ$-module equipped with the Mukai pairing; 
a {\it marked family} is a smooth family 
$\pi\colon\X\rightarrow S$ with topological fibre
$X$ and with a fixed isomorphism $\phi$ of the local system 
${\rm R}^*\pi_*\dZ$ with the constant local system $\tilde L$ 
on $S$ respecting pairings.

A {\it kernel} (derived correspondence) between smooth projective varieties 
$X_i$, $i=1,2$ is by definition an object $U\in\db{X_1\times X_2}$.
There is a composition product on kernels given for  
$U\in\db{X_1\times X_2}$ and $V\in\db{X_2\times X_3}$ by the standard formula
\[ U\circ V= {\bf R}p_{13*} \left({\bf L}p_{12}^*(U)\stackrel{{\bf L}}{\otimes} {\bf L}p_{23}^*(V)\right)\in \db{X_1\times X_3};
\]
here $p_{ij}\colon X_1\times X_2\times X_3\rightarrow X_i\times X_j$ 
are the projection maps. A kernel $U\in\db{X_1\times X_2}$ 
is {\it invertible}, if there is a kernel $V\in\db{X_2\times X_1}$ 
such that the products $U\circ V$ and $V\circ U$ are 
isomorphic in $\db{X_i\times X_i}$ to $\sheafO_{\Delta_{X_i}}$, the (complexes
consisting of) the structure sheaves of the diagonals. 

A kernel $U\in\db{X_1\times X_2}$ defines a functor
\[\Psi^U\colon \db{X_2}\rightarrow\db{X_1}
\]
by 
\[\Psi^U(-)= {\bf R}p_{1*}(U\stackrel{{\bf L}}{\otimes} p_2^*(-)), 
\]
If $U$ is invertible then $\Psi^U$ is a {\it Fourier--Mukai functor},
an equivalence of triangulated categories. 
An invertible kernel also induces an isomorphism
\[\psi^U\colon H^*(X_2, \dQ) \rightarrow H^*(X_1, \dQ)
\]
on cohomology, defined by 
\[ \psi^U(-)=p_{1*}(v(U)\cup p_2^*(-)). 
\]
$\Psi^U$ and $\psi^U$ are compatible via the Mukai map $v$. For 
Calabi--Yau manifolds of dimension at most four, $\psi^U$ is an isometry 
with respect to the Mukai pairing. 

Finally if $X_0$ is a complex manifold and 
$\alpha\in H^2(X_0, \sheafO_{X_0}^*)$, then 
there is a notion of an {\it $\alpha$-twisted sheaf} on $X_0$. This can either
be thought of as a sheaf living on the gerbe over $X_0$ defined by $\alpha$, 
or more explicitly as a collection
of sheaves on an open cover where the gluing
conditions are twisted by a Cech representative of $\alpha$. 

\section{Mirror symmetry and diffeomorphisms}
\label{mirrors}

The mirror symmetry story begins with the data $(Y, I_0, \omega_0)$, 
where $Y$ is a differentiable manifold, $I_0$ is a Calabi--Yau complex 
structure on $Y$, and $\omega_0$ is a K\"ahler\ form in this complex structure 
giving rise to a unique Ricci-flat metric\footnote{The $B$-field is discussed
in Remark~\ref{Bfield} below.}. 
The mirror of $(Y, I_0, \omega_0)$ is conjecturally found using the 
following procedure. 

\begin{procedure} \ \ \rm

\noindent {\bf Step 1} \ 
Choose a degeneration $\Y\rightarrow \Delta^*$ of $Y_0$
over a punctured multi-disc $\Delta^*$ with mid-point $P\in\Delta$, which is
a complex cusp (large complex structure limit point), 
an intersection of boundary divisors of the complex moduli space
with prescribed Hodge theoretic behaviour~\cite{morrison_compactifications}.  

\vspace{0.05in}
\noindent {\bf Step 2} \  
The choice of complex cusp $P$ induces a torus fibration 
on $Y_0$ via a degeneration process~\cite{gross_wilson2},~\cite{konts_soi}. 
Arrange this fibration so that its fibres are (special) Lagrangian with respect
to the complex structure $I_0$ and the symplectic form $\omega_0$. 

\vspace{0.05in}
\noindent {\bf Step 3} \ 
Dualizing this fibration (Strominger--Yau--Zaslow~\cite{syz}) should 
yield, after an appropriate compactification, a 
compact differentiable manifold $X$, together with an isomorphism
\[
{\rm mir}_P\colon H^*(Y, \dQ) \rightarrow H^*(X,\dQ), 
\]
relating algebraic cohomology to 
transcendental cohomology, and compatible with filtrations coming 
from the Leray spectral sequence of the torus fibration on the 
transcendental part and degree on the algebraic part. 

\vspace{0.05in}
\noindent {\bf Step 4} \ The symplectic form 
$\omega_0$ and complex structure $I_0$ on $Y$ induce 
a Calabi--Yau structure $(X, I'_0)$ and a symplectic form
$\omega'_0$ on the manifold~$X$. 

\vspace{0.05in}
\noindent {\bf Step 5} \  
According to Kontsevich' homological mirror symmetry 
conjecture~\cite{kontsevich} (cf. also Manin~\cite{manin}), 
the statement that $(Y, I_0, \omega_0)$ and $(X, I'_0, \omega'_0)$ 
are mirror symmetric is expressed by a equivalence of categories
determined by the torus fibration, 
\[\mir_P\colon {D^b{\rm Fuk}}(Y, \omega_0) \longrightarrow \db{X_0}.  
\]
On the left hand side, ${D^b{\rm Fuk}}(Y, \omega_0)$ is the derived 
Fukaya category, which is (conjecturally) constructed from  
Lagrangian submanifolds and flat bundles on them. 
This equivalence is compatible with the cohomology isomorphism 
${\rm mir}_P$ of Step 3. 
\label{mirrorproc}
\end{procedure}

My main interest in this paper lies in an understanding of the action 
of the diffeomorphism group $\diff^+(Y)$ of $Y$ on the above data. 
A diffeomorphism $\gamma\in\diff^+(Y)$ maps a triple $(Y, I_0, \omega_0)$ 
to a new triple $(Y, \gamma^*I_0, \gamma^*\omega_0)$ and in particular it 
defines a symplectomorphism $(Y, \gamma^*\omega_0)\cong (Y, \omega_0)$. 
Lifting $\gamma$ to a graded symplectomorphism~\cite{seidel_gr} 
induces an equivalence
\[ \tilde\gamma\colon {D^b{\rm Fuk}}(Y, \gamma^*\omega_0) \stackrel{\sim}{\longrightarrow}  {D^b{\rm Fuk}}(Y, \omega_0) \ {\rm mod}\,[1]\]
i.e. well-defined up to translation; this comes from the fact that the
lifting of $\gamma$ to a graded symplectomorphism is only well-defined
up to shift. 

Perform Procedure~\ref{mirrorproc} on the pair of Calabi--Yau triples 
$(Y, I_0, \omega_0)$ and $(Y, \gamma^*I_0, \gamma^*\omega_0)$ 
with respect to the same complex cusp~$P$ to obtain mirrors  
$(X, I'_0, \omega'_0)$ and $(X, I''_0, \omega''_0)$. Then there is a 
diagram of categorical equivalences up to translation, where the bottom 
arrow is defined by the others: 

\[\begin{array}{ccc}
{D^b{\rm Fuk}}(Y, \gamma^*\omega_0) & \stackrel{\tilde\gamma}{\longrightarrow} & {D^b{\rm Fuk}}(Y, \omega_0) \\
\mapdownleft{\!\!\!\!\!\!\!\!\!\mir_P} && \mapdownright{\mir_P} \\
\db{X, I''_0} & \stackrel{\Psi_0}{\longrightarrow} & \db{X, I'_0}
\end{array}
\]
Further, the categorical equivalence $\Psi_0$ induces 
a cohomology isomorphism 
$\psi_0$ which must be compatible via the isomorphism ${\rm mir}_P$ with the action 
of $\gamma^*$ on $H^*(Y,\dQ)$. 

\vspace{0.05in}

To push this further, assume that $(I_t, \omega_t)$ is a family 
of complex structures and compatible symplectic forms on $Y$ 
deforming $(I_0, \omega_0)$. 
The diffeomorphism $\gamma$ induces a family of equivalences
\[ \tilde\gamma\colon {D^b{\rm Fuk}}(Y, \gamma^*\omega_t) \stackrel{\sim}{\longrightarrow}  {D^b{\rm Fuk}}(Y, \omega_t) \ {\rm mod}\,[1].\]
The mirrors of $(Y, I_t, \omega_t)$ and $(Y, \gamma^*I_t, \gamma^*\omega_t)$ 
with respect to the cusp $P$ are deformations
$(X, I_t', \omega_t')$ and $(X, I''_t, \omega_t'')$ 
of the original $(X, I_0', \omega_t')$ and $(X, I''_0, \omega_t'')$. 
Consequently, there are induced equivalences
\[
\Psi_t\colon \db{X,I_t''}  \stackrel{\sim}{\longrightarrow}  \db{X,I_t'} \ {\rm mod}\,[1].
\]

\begin{slogan} Diffeomorphisms of the manifold $Y$ induce families of 
equivalences of derived categories up to translation over 
the complex moduli space of the mirror manifold $X$.
\label{mainslogan}
\end{slogan}

In the next Section, this slogan will be translated into a precise conjecture. 
The rest of the paper is providing evidence and examples. Before that 
however, there are two important points to clear up. 

\begin{remark} \rm  In order to be able to compare cohomology actions, 
it is necessary to work with marked moduli spaces of 
Calabi--Yau varieties in the procedure described above. 
This implies that the moduli spaces arising will typically
be non-connected; this issue already arises for K3 surfaces. The
procedure outlined above tacitly assumed that for $\gamma\in\diff^+(Y)$, 
the complex structures $I_0$ and $\gamma^*I_0$ will be in the same connected
component of the (marked) moduli space of complex structures on $Y$, and thus
the same complex cusp $P$ may be used to form the mirror. In general, 
several components of the marked moduli space will contain complex structures
related by diffeomorphisms. The way out is to pick a full set of components of
the complex structure moduli space of $Y$ containing structures related by 
diffeomorphisms, and fix a complex cusp $P_i$ in each one of them 
in a compatible way so that the topological mirror $X$ and the cohomology map
${\rm mir}_P$ are the same for all cusps. 
Then the procedure described above is accurate. I~tacitly assume
this extension of the setup everywhere below. 
\label{cusps} 
\end{remark}

\begin{remark} \rm The above discussion ignores one crucial
piece of data present in mirror symmetry, namely the $B$-field. In this remark
I want to argue that, possibly under an extra assumption in the case of 
threefolds, I can consistently restrict to the case of vanishing $B$-field
on the holomorphic side of mirror symmetry. The generalization of the
present ideas to the case of arbitrary $B$-fields is left for future work. 

Physical mirror symmetry deals with quadruples $(Y, I_0, \omega_0, B_0)$; 
here $B_0\in H^2(Y, \dR/\dZ)$ is the $B$-field on the Calabi--Yau manifold
$(Y, I_0)$. The mirror of $(Y, I_0, \omega_0, B_0)$ is a quadruple 
$(X, I'_0, \omega'_0, B'_0)$. According to a physics proposal originally 
formulated in the context of K-theory, for nonzero $B$-fields Kontsevich' 
homological mirror symmetry conjecture should take a form of an 
equivalence of categories  
\[ \mir_P\colon {D^b{\rm Fuk}}(Y, \omega_0, B_0) \longrightarrow \db{X, I'_0, B'_0}; \]
see for example Kapustin--Orlov~\cite{kapustin_orlov}. The category
${D^b{\rm Fuk}}(Y, \omega_0, B_0)$ on the left hand side  
should be a derived Fukaya-type category where the bundles
on Lagrangian submanifolds are equipped with connections whose 
curvature is given by the restriction of $B$. For the current picture
this is not a serious issue as diffeomorphisms should still 
give equivalences between these generalized derived Fukaya categories. 

More importantly, the category on the right hand 
side should be a version of the 
derived category of coherent sheaves on $(X, I'_0)$. 
More precisely, consider the natural map 
\[ \delta\colon H^2(X_0, \dR/\dZ)\rightarrow H^2(X_0, \sheafO_{X_0}^*). 
\]
Suppose first that the class $B$ has torsion image in 
$H^2(X_0,\sheafO_{X_0}^*)$. 
In this case
 $\db{X_0, B'_0}$ should be the derived category of $\delta(B'_0)$-twisted
coherent sheaves $\db{X_0, \delta(B'_0)}$. If $\delta(B)$ is non-torsion, 
then $\db{X_0, B'_0}$ should be a suitable subcategory of a category of 
sheaves over a gerbe; for more discussion
see~\cite[Remark 2.6]{kapustin_orlov}. 

Note that for elliptic curves, $H^2(X,\sheafO_X^*)=0$ from the exponential 
sequence, so the derived category is not affected. However in higher 
dimensions, there are nontrivial twists. 

Consider the case of K3 surfaces. The induced action of $\diff^+(Y)$ 
on the space parameterizing quad\-ru\-ples $(X, I', \omega', B')$
is compatible via the cohomology isomorphism ${\rm mir}_P$ with 
the action on the rational cohomology $H^*(Y, \dQ)$. 
Moreover, by Proposition~\ref{K3eta}, the cohomology isomorphism ${\rm mir}_P$ is 
defined over $\dZ$. So the induced action of $\diff^+(Y)$ maps a zero 
$B'$-field to a zero $B'$-field. 

In the case of Calabi--Yau threefolds, the isomorphism ${\rm mir}_P$ is not 
expected to be defined over $\dZ$ but it is still defined over $\dQ$ so 
it preserves torsion $B'$-fields. On the other hand, in this case the 
torsion subgroup of $H^2(X,\sheafO_X^*)$ is isomorphic to the torsion 
subgroup of $H^3(X,\dZ)$. Assuming $H^3(X,\dZ)_{\rm tors}=0$, 
there are no possible torsion twists at all. As there is no known 
example of a smooth Calabi--Yau threefold with torsion in $H^3$, this 
restriction does not appear to be very serious. 

I conclude therefore that for elliptic curves, K3 surfaces and 
Calabi--Yau threefolds under the assumption $H^3(X,\dZ)_{\rm tors}=0$, 
it is legitimate to set $B'=0$, as was effectively done in the above
discussion. 
\label{Bfield}
\end{remark}

\section{The categorical mapping group} 

To translate Slogan~\ref{mainslogan} into a conjecture, I
formalize the notion of a family 
of equivalences for complex structures on $X$. 
Assume that $\pi\colon\X\rightarrow \D$ is a family of complex projective 
varieties over a smooth complex base $\D$. Consider triples $(\phi,\alpha,U)$
where  
\begin{itemize}
\item $S$ is an open subset of $\D$, the complement of a countable number of 
closed analytic submanifolds, 
\item $\phi\colon S\rightarrow \D$ is an analytic injection (not necessarily the identity), giving rise to the fibre product diagram 
\[\begin{array}{ccc}
\X_\phi & \longrightarrow & \X_S \\
\Big\downarrow && \mapdownright{\pi} \\
\X_S & \stackrel{\phi\circ\pi}{\longrightarrow} &  S
\end{array}\]
where $\pi\colon\X_S\rightarrow S$ is the restriction of the original family
to $S$; 
\item $\alpha\in H^2(S, \sheafO_S^*)$ is a twisting class on the base;  
\item $U$ is an object of the bounded derived category of quasi-coherent
$p^*\alpha$-twisted sheaves on $\X_\phi$, whose derived restriction to 
the fibres of $p$ is isomorphic to a bounded complex of coherent sheaves
on $X_s\times X_{\phi(s)}$; here $p\colon\X_\phi\rightarrow S$ 
is the natural map.
\end{itemize}
Let $M(\X, \D)$ be the set of triples under obvious identifications, namely
extension of $\phi$, and tensoring by a line bundle 
pulled back from the base $S$. 
$M(\X, \D)$ can be given a monoid structure by a standard procedure
generalizing the composition structure on kernels; the maps $\phi$
simply compose under the multiplication rule. 
The multiplication has a two-sided unit $({\rm id}_S, 0, \sheafO_{\Delta_\X})$ 
where $\Delta_\X\subset \X\times_S \X$ is the relative diagonal. 

Suppose further that $\X\rightarrow \D$ is a marked family
with topological fibre $X$, and let 
$\tilde L_\dQ=H^*(X,\dQ)$ be the rational cohomology of the fibres.

\begin{definition} The {\rm categorical mapping group} 
$G(\X, \D)$ of the marked family
$\pi\colon\X\rightarrow \D$ is the group of invertible elements $g\!=\!(\phi,\alpha,U)$
of the monoid $M(\X, \D)$ satisfying the following cohomological 
condition\footnote{Note that the cohomology condition is vacuous 
if $\D$ is connected.}:
there should exist an isometry $\psi_g\in\aut(\tilde L_\dQ)$ 
and a commutative diagram of local systems on $S$  
\[ \begin{array}{ccc} 
R^*(\phi\circ \pi)_*\dQ_\X & \longrightarrow & R^*\pi_*\dQ_\X\\
\downarrow && \downarrow \\
\tilde L_\dQ & \stackrel{\psi_g}{\longrightarrow} &\tilde L_\dQ. 
\end{array} 
\]
Here the vertical arrows come from the
markings, and the top horizontal arrow is the
isomorphism of local systems induced by the invertible relative kernel $U$. 

The image $G^{coh}(\X, \D)$ of the resulting homomorphism 
\[ \tilde \psi\colon G(\X, \D)\rightarrow \aut(\tilde L_\dQ). 
\]
is the {\rm cohomological mapping group} of
the family $\X\rightarrow \D$.
\end{definition}

For an element $g=(\phi,\alpha,U)\in G(\X, \D)$,
derived restriction of $U$ to the fibre over $s\in S$
gives a well-defined untwisted object $U_s$ in $\db{X_s\times X_{\phi(s)}}$. 
An easy base change argument shows that as $(\phi,\alpha,U)$ is invertible, 
$U_s$ is also invertible, and hence there is a Fourier--Mukai transform
\[ \Psi^{U_s}\colon\db{X_{\phi(s)}} \rightarrow \db{X_{s}}
\]
for every $s\in S$, with a well-defined cohomology action $\psi_g$ independent
of $s\in S$. 

The following main conjecture summarizes the discussion of the preceding 
sections.  

\begin{conjecture}
\label{mainconjecture}
Let $(\Y\rightarrow T,P)$ be a family of Calabi--Yau manifolds together with a 
choice of complex cusp $P\in\partial T$. Let $X$ be the topological 
mirror of $Y$ with respect to $P$; assume that $H^3(X,\dZ)_{\rm tors}=0$.  
Let $\X\rightarrow \M_X$ be the marked
Calabi--Yau moduli space (Teichm\"uller space) of $X$. Then there exists 
a homomorphism $\xi$ fitting into a commutative diagram
\[ \begin{array}{ccc}
\diff^+(Y)  & \stackrel{\xi}{\longrightarrow}& G(\X, \M_X)/{\rm (translations)} \\
\mapdownleft{\beta}&&\mapdownright{\psi}\\
\aut(H^*(Y,\dQ)) & \longrightarrow & \aut(H^*(X,\dQ))/(\pm 1),
\end{array}
\]
where 
\begin{itemize}\item translations mean relative translations over connected components on
the base $\M_X$; moreover, 
\item$\beta$ is the natural action,  
\item$\psi$ is the map associating to an element $g$ of the categorical 
mapping group its cohomology action $\psi_g$ and 
\item the bottom horizontal arrow 
is the canonical map 
induced by the isomorphism ${\rm mir}_P$ of Procedure~\ref{mirrorproc}.
\end{itemize}

The diagram gives rise to an isomorphism of groups 
\[ \bar\xi\colon \langle\diff^{\rm coh}(Y), (-1)\rangle/(\pm 1) {\longrightarrow} G^{\rm coh, filtr}(\X, \M_X)/(\pm 1)
\]
where 
\begin{itemize}\item $\diff^{\rm coh}(Y)$ is the image of $\diff^+(Y)$ in 
$\aut(H^*(Y,\dQ))$, and
\item $G^{\rm coh, filtr}(\X, \M_X)$ is the subgroup of 
$G^{\rm coh}(\X, \M_X)$ consisting of elements whose associated cohomology 
action preserves the Leray filtration on the transcendental cohomology of~$X$.
\end{itemize}
\end{conjecture} 

\begin{remark} \rm 
As discussed before, a diffeomorphism only acts on the derived Fukaya 
category up to translation. Hence in going from diffeomorphisms to their
cohomology action, I have to take account of the cohomology action of
translation $[1]$ on the derived
Fukaya category which is simply multiplication by 
$(-1)$. This is reflected on the left hand side of the cohomological form of
the conjecture. On the right hand side, the cohomology actions of relevant
families of Fourier--Mukai transforms must preserve the Leray filtration on 
the transcendental part of the cohomology of $X$; 
this is mirror to the statement that 
diffeomorphisms preserve the degree filtration on 
the algebraic cohomology of $Y$.
\end{remark} 

\section{Elliptic curves}

The diffeomorphism group $\diff^+(T^2)$ of the two-torus
acts via the standard 
$SL(2, \dZ)$-action on $H^1(T^2, \dZ)$ and trivially on $H^2(T^2, \dR)$.  
Thus, recalling the discussion of Section~\ref{mirrors}, $\diff^+(T^2)$
should be realized by self-equi\-va\-len\-ces in the mirror family. 

The elliptic curve is self-mirror. Choose cohomology classes 
dual to a collapsing circle fibre and a section of a 
circle fibration corresponding to the decomposition $T^2=(\dR/\dZ)^2$
to fix an isomorphism 
\begin{equation}
{\rm mir}_P\colon H^*(T^2,\dZ) \rightarrow H^*(T^2, \dZ)
\label{elliptic_iso}
\end{equation}
interchanging even and odd cohomology. Let 
\[\X=\dC\times{\mathcal H}/[(z, \tau)\sim (z+n+ m\tau,\tau)], \ \   m,n \in \dZ
\]
be the universal family of marked elliptic curves together with the map 
\[\pi\colon \X\rightarrow {\mathcal H}\]
to the marked moduli space, the upper half plane ${\mathcal H}$. 

\begin{theorem} There exists a group homomorphism 
\[ \xi\colon \diff^+(T^2)\rightarrow G(\X,{\mathcal H})/{\rm (translations)}
\]
compatible with the isomorphism ${\rm mir}_P$ and descending to an isomorphism
\[\bar\xi\colon PSL(2,\dZ)\stackrel{\sim}\rightarrow G^{\rm coh, filtr}(\X,{\mathcal H})/(\pm 1).
\]
\label{ellipticthm}
\end{theorem}
{\vspace{.05in}\noindent {\sc Sketch Proof} \hspace{.05in}}
One only has to relativize the discussion of~\cite[Section 3d]{st}; the 
$PSL(2,\dZ)$-action on the derived category of an elliptic curve was
of course already known to Mukai~\cite{mukai_abelian}. 
There are no nontrivial twists over the base ${\mathcal H}$ 
and for all elements $(\phi, U)\in G(\X, {\mathcal H})$ in the image of 
$\xi$, the map $\phi$ is the identity. 
\ethrm 

\section{K3 surfaces}

\subsection{Lattice polarized mirror symmetry}
\label{K3s}

Let $M$ be an even non-degenerate sublattice of signature $(1,t)$ of 
the K3 lattice $L=H^2(X,\dZ)$. A {\it marked ample $M$-polarized 
K3 surface} is a marked $K3$ surface $(Y_t,\psi)$, which satisfies 
$\psi^{-1}(M)\subset \pic(Y)$ and moreover $\psi^{-1}({\mathcal C})$ contains
an ample class for ${\mathcal C}$ a chosen chamber of $M_\dR$; 
for the precise definition see~\cite{dolgachev}. 
For a pair $(Y_t,\psi)$, Hodge decomposition on the second cohomology 
gives the period point in 
\[ \D_M=\left\{ z\in \dP(M^\perp_\dC) : z^2 =0, z\cdot \bar z > 0  \right\}\setminus \bigcup_{\delta\in\Delta(M^\perp)} \langle\delta\rangle^\perp \]
where $\Delta(M^\perp)$ is the set of vectors in $M^\perp$ of length $-2$. 
The period domain $\D_M$ has two connected components. 

\begin{proposition} {\rm (\cite[Corollary 3.2]{dolgachev})}
The period map realizes $\D_M$ as the marked moduli space of 
marked ample $M$-polarized K3 surfaces; in particular, there exists 
a universal marked family $\pi\colon\Y\rightarrow \D_M$. 
\end{proposition}

The mirror construction works best for a sublattice $M$ of $L$ with
$M^\perp$ containing hyperbolic planes; I assume this from now on. The
choice of a pair of vectors $\{f, f'\}\subset M^\perp$ spanning a hyperbolic 
plane amounts to choosing a complex cusp\footnote{The choice in fact 
specifies a complex cusp in both 
components of $\D_M$; c.f. Remark~\ref{cusps}.} 
$P\in \partial \D_M$. The choice of $\{f, f'\}$ in $M^\perp$
also induces a splitting
$M^\perp = \check M \perp \langle f, f'\rangle$. This defines 
a new sublattice $\check M$ of $L$  
of signature $(1, 20-t)$ and gives rise to a family 
$\X\rightarrow \D_{\check M}$. This is the mirror family of 
$\Y\rightarrow \D_M$. The mirror construction, under the assumption on $M$, 
is an involution: starting with the pair $\{f, f'\}$ thought of as defining
a complex cusp in the boundary of $\D_{\check M}$, one recovers the original 
$M$-polarized family $\Y\rightarrow \D_M$.
Monodromy and filtration considerations show

\begin{proposition}\label{K3eta} The cohomology isometry
\[ {\rm mir}_P\colon H^*(Y, \dQ)\rightarrow H^*(X, \dQ), 
\]
where $Y$, $X$ are thought of as fibres in $M$-, respectively 
$\check M$-polarized families, is given by 
\[\begin{array}{c}
{\rm mir}_P(h_0) = f',\ \  {\rm mir}_P(h_4) = -f, \ \ {\rm mir}_P(f)  =  -h_4,\ \  {\rm mir}_P(f') =  h_0, \\
{\rm mir}_P(m)  =  m \mbox{ for } m\in \langle h_i, f_i \rangle^\perp.
\end{array}\] 
In particular ${\rm mir}_P$ is defined over $\dZ$. 
\end{proposition} 

In the framework of Procedure~\ref{mirrorproc}, I want to start with a general 
ample $M$-polarized K3 $(Y, I_0, \omega_0)$ and I want to use the degeneration 
to the complex cusp $P\in\partial\D_M$
defined by $\{f, f'\}$ in the $M$-polarized family 
$\Y\rightarrow \D_M$. The cohomology class of the 
the K\"ahler\ form has to be orthogonal to the period, so 
$\omega_0$ is restricted to live in the subspace 
$M_\dR\subset  H^2(Y,\dR)$. This condition has to be preserved throughout
the whole procedure, hence the set of diffeomorphisms has 
to be restricted as well. Call $\gamma\in\diff^+(Y)$ {\it $M$-allowable}, 
if $\gamma^*\in O(L)$ satisfies $\gamma^*(M)=M$. Let $\diff_M^+(Y)$ be the 
group of $M$-allowable diffeomorphisms of $Y$ and 
$\diff_M^{\rm coh}(Y)$ the image of $\diff_M^+(Y)$ in the isometry 
group of $\tilde L$. Note that monodromy transformations around various 
rational boundary components of the marked moduli space $\D_M$ are certainly 
$M$-allowable. 

\subsection{Action on cohomology} 

Before I proceed further, let me recall

\begin{theorem} {\rm (Donaldson~\cite{donaldson})}
Let $\gamma\in\diff^+(Y)$ be an orientation-preserving diffeomorphism 
of a K3 surface. 
Then the induced action on cohomology preserves the orientation of 
positive definite three-planes in $H^2(Y,\dR)$. 
\label{donthm}
\end{theorem}

Comparing this theorem and the cohomology isomorphism ${\rm mir}_P$ given above
leads to an important observation. 
Fix once and for all an orientation for the positive part 
of $H^0(X)\oplus H^4(X)$. Then for an algebraic K3 surface $X_0$
this, together with the real and imaginary parts of the 
period of $X_0$ and an ample class, gives an orientation to 
positive definite four-planes in the inner product space
$H^*(X_0, \dR)$, that I call for want of a better expression the
{\it canonical orientation}. If the diffeomorphism group maps, at least on
the level of cohomology, surjectively onto the cohomological mapping 
group, then Donaldson's theorem must have an analogue for Fourier--Mukai 
functors. I formulate this as\footnote{A weaker form of this conjecture
was proposed independently by Markman~\cite{markman}.}

\begin{conjecture}\label{orient} Let $\Psi\colon \db{Z_0} \rightarrow \db{X_0}$
be a Fourier--Mukai equivalence between derived categories of
smooth projective K3 surfaces. 
Then the induced cohomology action $\psi\colon H^*(Z_0,\dZ) \rightarrow H^*(X_0,\dZ)$ 
preserves the canonical orientation of positive definite four-planes 
in cohomology. 
\end{conjecture}

%It is tempting to speculate that if true,  this conjecture could be proved using an invariant `counting' sheaves, analogous to the Donaldson invariant in the proof of Donaldson's theorem. For Calabi--Yau threefolds, such an invariant isconstructed by Donaldson and Thomas~\cite{dt},~\cite{thomas}. However, I do not know how to go about doing this. 

\subsection{Cohomological realization of the main conjecture}

All the pieces are together; the following theorem confirms that
Conjecture~\ref{mainconjecture} can be realized on a cohomological level
for mirror families of polarized K3 surfaces. 

\begin{theorem} Assume Conjecture~\ref{orient}. 
Let $P$ be a complex cusp in the boundary $\partial \D_M$
of the marked moduli space of the family of ample $M$-polarized K3 surfaces. 
Let $\X\rightarrow \D_{\check M}$ be the mirror family. Then 
there exists an isomorphism
\[ \bar\xi\colon \left\langle\diff^{\rm coh}_M(Y), (-1)\right\rangle\Big/(\pm 1) \stackrel{\sim}\longrightarrow G^{\rm coh, filtr}(\X, \D_{\check M})/(\pm 1)
\]
compatible with the isometry ${\rm mir}_P$. 
\label{K3_main_theorem}
\end{theorem}
{\vspace{.05in}\noindent {\sc Sketch Proof} \hspace{.05in}}
Let
\[\sigma\in{\rm mir}_P\, \left\langle\diff^{\rm coh}_M(Y), (-1)\right\rangle\, {\rm mir}_P^{-1}\subset O(\tilde L)\]
be an isometry of the Mukai lattice $\tilde L$. The task is to find an element
$(\phi, \alpha, U)\in G(\X, \D_{\check M})$ with cohomology action $\sigma$.  

Proposition~\ref{K3eta} shows that 
$\sigma$ fixes $\check M^\perp = M\perp \langle f, f'\rangle$. 
Therefore the map $\sigma^{-1}$ induces an automorphism 
\[\phi\colon \D_{\check M} \rightarrow \D_{\check M}\] 
mapping $[z]\in\D_{\check M}$ to $[\sigma^{-1}_\dC(z)]$. 
The point about this map is that 
for $s\in\D_{\check M}$, $\sigma$ induces a Hodge isometry
from $H^*(X_{\phi(s)},\dZ)$ to $H^*(X_s,\dZ)$ that
by Donaldson's Theorem~\ref{donthm} preserves the canonical orientation.

Work of Mukai~\cite{mukai} and global Torelli
identifies $X_{\phi(s)}$ with a moduli
space of sheaves on $X_t$ and this can be done in a relative way 
over an open subset $S\subset\D_{\check M}$ in the family 
$\X\rightarrow\D_{\check M}$. Some further arguments show that there exists
a (suitably twisted) universal sheaf $U$ which defines an element 
of $G(\X,\D_{\check M})$ with cohomology action $\sigma$.  
\ethrm 

\subsection{Toric families}
\label{toricK3}

I want to give some examples of the correspondence between 
diffeomorphisms and families of Fourier--Mukai transformations in toric
families. Let $\Delta\subset\dZ^3\otimes \dR$ be a three-dimensional reflexive
lattice polyhedron~\cite{batyrev}. It gives rise to a toric
threefold $\dP_\Delta$, which has a crepant toric resolution 
$\tilde\dP_\Delta\rightarrow\dP_\Delta$. 

Let $\X_\Delta\rightarrow\F_\Delta$ be the family  
of anticanonical hypersurfaces 
$X_t\subset\tilde\dP_\Delta$. For $t\in\F_\Delta$, 
$X_t$ is a smooth K3 surface and 
there is a restriction map $\pic(\tilde\dP_\Delta)\rightarrow \pic(X_t)$
with image $M_\Delta\subset\pic(X_t)$, the space of toric
divisors on $X_t$. The equality $M_\Delta=\pic(X_t)$ for general $t$
holds only under an extra condition, 
namely if for all one-dimensional faces $\Gamma$ of $\Delta$, neither 
$\Gamma$ nor its dual $\Gamma^*\subset\Delta^*$ contain extra lattice 
points in their interior. If this condition fails, and it frequently does, 
then the Picard lattice of $X_t$ is not spanned by toric divisors. 

Associated to the lattice $M_\Delta$ is a family of ample $M_\Delta$-polarized
surfaces $\X\rightarrow\D_{M_\Delta}$; 
$\X_\Delta\rightarrow\F_\Delta$ is a subfamily 
of this family. However, if the above condition fails, then the general ample 
$M_\Delta$-polarized K3 is not toric and 
$\X_\Delta\rightarrow\F_\Delta$ is a proper 
subfamily of $\X\rightarrow \D_{M_\Delta}$.

Let $\Delta^*$ be the dual polyhedron of $\Delta$, and let
$M_{\Delta^*}$ be the lattice of toric divisors of the Batyrev
mirror family $\Y_{\Delta^*}\rightarrow\F_{\Delta^*}$. 
According to a conjecture of Dolgachev~\cite{dolgachev}, there should exist
a decomposition $M_\Delta^\perp=\langle f, f'\rangle\perp\check M_\Delta$
as before, together with a primitive embedding 
\[ i\colon M_{\Delta^*}\hookrightarrow \check M_\Delta. 
\]
This conjecture is supported by numerical evidence and 
can be checked in several concrete cases~\cite[Section 8]{dolgachev}; I 
assume that $\Delta$ is chosen so that the embedding exists.  
The embedding $i$ realizes inside the period domain
\[\D_{M_\Delta}\subset\dP(M_\Delta^\perp)_\dC=\dP(\check M_\Delta\perp \langle f, f'\rangle)_\dC\]
the toric deformation space $\F_\Delta$ as an explicit subdomain  
\[\F_\Delta=\D_{M_\Delta}\cap \dP(M_{\Delta^*}\perp \langle f, f'\rangle)_\dC\]
(cf. \cite[Section 4.3]{kobayashi}). 

The decomposition $M_\Delta^\perp=\langle f, f'\rangle\perp \check M_\Delta$
gives the family $\Y\rightarrow\D_{\check M}$ of marked 
ample $\check M_\Delta$-polarized surfaces as the mirror family of
$\X\rightarrow\D_{M_\Delta}$. 

\vspace{0.2in}

Turning to the specific issue of diffeomorphisms, 
Theorem~\ref{K3_main_theorem} realizes, on the level of cohomology, 
$\check M_\Delta$-allowable diffeomorphisms 
on the family $\X\rightarrow\D_{M_\Delta}$ by Fourier--Mukai 
functors\footnote{A word of warning: the notation used here is related to that used in Theorem~\ref{K3_main_theorem} by $M=\check M_\Delta$ and consequently $\check M=M_\Delta$.}. 
The first examples are $\check M_\Delta$-allowable diffeomorphisms arising as 
monodromy transformations around boundary components meeting at the 
large complex structure limit point in the boundary of $\D_{\check M_\Delta}$. 

\begin{lemma} {\rm (\cite[6.2]{dolgachev})}
The monodromy transformations in boundary components meeting at the complex 
cusp $P\subset\partial\D_{\check M_\Delta}$ 
generate a group isomorphic to $M_\Delta$. 
The transformation corresponding to $d\in M_\Delta$ acts on
the K3 lattice $L=H^2(Y,\dZ)$ by 
\def\arraystretch{1.25}
\[\begin{array}{lcl} T_d(f) & = & f \\
T_d(f') & = & f' + d - \frac{d^2}{2} f \\
T_d(m) & = & m - (d\cdot m) f \mbox{ for } m\in \langle f, f' \rangle^\perp.
\end{array}\] 
\end{lemma} 

This action fixes the lattice 
$\check M_{\Delta}\perp\langle h_0, h_4\rangle$ pointwise, 
so a glance at the definition the 
map $\phi$ in the sketch proof of Theorem~\ref{K3_main_theorem} 
shows that $\phi$ acts trivially on $\D_{M_\Delta}$. Indeed, 
by common wisdom, the monodromy transformation $T_d$ 
is mirrored by tensoring by the line bundle $\sheafO_{X_t}(d)$ 
for $d\in M_\Delta\subset\pic(X_t)$ and $t\in\D_{M_\Delta}$. 
In the relative framework, this is simply tensoring by a relative 
line bundle (twisted over the base if necessary). 

However, as one moves away from these boundary components, 
diffeomorphisms with more complicated cohomology action crop up. 
In particular, there are  diffeomorphisms 
$\gamma\in\diff^+_{\check M_\Delta}(Y)$
whose associated cohomology action $\gamma^*\subset O(L)$ fixes every element 
of $M_{\Delta^*}$ but does not fix the whole of 
$\check M_{\Delta}$ pointwise. Hence these diffeomorphisms have mirrors with a 
nontrivial map $\phi\colon\D_{M_\Delta}\rightarrow \D_{M_\Delta}$. 
This nontrivial map has a fixed locus which includes 
the toric locus $\F_\Delta\subset\D_{M_\Delta}$. The 
interpretation is that these diffeomorphisms are mirrored by 
families of Fourier--Mukai transforms acting as self-equivalences 
over the toric
locus $\F_\Delta$, but acting by equivalences between derived categories 
of different K3s away from the toric locus. Compare Example~\ref{typeiiiex}. 

\section{Threefolds} 

\subsection{Examples of families of self-equivalences} 

\begin{example}\rm 
Suppose that $\X\rightarrow \M_X$ is a marked family of Calabi--Yau 
threefolds. As in the previous section, some immediate examples of elements of
$G(\X, \M_X)$ come from line bundles. Let $d\in H^2(X, \dZ)$; in the
threefold case $d$ is automatically the first Chern class 
of a line bundle $\sheafO_{X_s}(d)\in\pic(X_s)$
on the fibres $X_s$ for $s\in\M_X$. There is a relative 
line bundle $\sheafO_\X(d)$, possibly twisted over the base, which when 
thought of as a sheaf on the diagonal $\Delta_\X\subset\X\times_{\M_X}\X$ 
gives an invertible relative kernel and so an element of $G(\X, \M_X)$ with 
$\phi$ the identity on $\M_X$. The fibrewise Fourier--Mukai transform
is tensoring by $\sheafO_{X_s}(d)$ in the derived
category of $X_s$. These elements of 
$G(\X, \M_X)$ should again mirror monodromy transformations around components
of the boundary of $\M_X$ meeting at the large structure limit point 
$P'\in\delta \M_X$. 
\end{example}

\begin{example}\rm\label{typei} A more interesting 
case where families of Fourier--Mukai self-equivalences arise 
is that of families of twist functors coming from contractions. 
Consider a Calabi--Yau threefold $X_0$
together with a flopping (or Type I) contraction $X_0\rightarrow \bar X_0$, a
birational contraction with exceptional locus of dimension one. 
There is a diagram 
\[\begin{array}{ccccc}{{\mathcal C}} & \rightarrow & \X & \rightarrow & \bar X\\
& \searrow & \downarrow & \swarrow\\
&& S
\end{array}\]
where $S\subset \M_X$ is a dense open subset in a 
connected component of the marked Calabi--Yau moduli space of $X_0$,
and ${{\mathcal C}}$ is flat over 
$S$ intersecting every fibre in a chain of rational 
curves. The structure sheaf $\sheafO_{{\mathcal C}}$
restricts to the fibre $X_s$ as $\sheafO_{{{\mathcal C}}_s}$,  
a spherical sheaf in the sense of Seidel--Thomas~\cite{st}. This sheaf
gives rise to a twist functor, a Fourier--Mukai self-equivalence 
of $\db{X_s}$. There is an invertible relative kernel $U$, 
giving rise to an element $(\phi, 0, U)$ of
$G(\X, \M_X)$ with $\phi$ the identical embedding of $S$ in $\M_X$ 
and restricting to the twisting kernel on $X_s\times X_s$ for $s\in S$. 
\end{example} 

\begin{example}\rm \label{typeii}
Suppose that a Calabi--Yau threefold\ $X_0$ has a Type II contraction, 
a birational morphism $X_0\rightarrow \bar X_0$ contracting an irreducible
divisor $E_0$ to a point. $E_0$ is possibly non-normal but in any case
it is a del Pezzo surface i.e. $\omega_{E_0}$ is ample. In this case, 
there is a diagram
\[\begin{array}{ccccc}{\mathcal E} & \rightarrow & \X & \rightarrow & \bar X\\
& \searrow & \downarrow & \swarrow\\
&& S
\end{array}\]
where $S\subset\M_X$ is dense in a connected component 
of the marked moduli space and the exceptional locus ${\mathcal E}$ is a flat 
family of surfaces over $S$ intersecting $X_t$ in a del Pezzo $E_s$. 

The structure sheaf $\sheafO_{E_s}$ of a del Pezzo surface is exceptional 
in the sense that there are no higher Ext's, 
so by~\cite[Proposition 3.15]{st}, 
its pushforward to $X_s$ is spherical. It defines a
twist functor, a self-equivalence of the derived category of 
$X_s$. The construction works again over $S$; there is an invertible 
relative kernel defined by the family of surfaces ${\mathcal E}$ over $S$ 
which gives an element of $G(\X, \M_X)$ with $\phi$ still being the identity. 
\end{example} 

The elements of $G(\X, \M_X)$ given in Examples~\ref{typei}--\ref{typeii}
should mirror monodromies around various 
other components of the boundary of moduli as discussed by Horja~\cite{horja}
and more recently by Aspinwall~\cite{aspinwall}. In particular, 
some of these transformations should mirror Dehn twists in Lagrangian
three-spheres in the mirror, which are known to arise in degenerations to 
threefolds with double points; this was one of the starting points for
Seidel and Thomas~\cite{st}. 

\subsection{An example with a nontrivial action on moduli}

In this section, I give an example of an element of the categorical mapping 
group with a nontrivial map $\phi$ on the base of moduli. The example comes
from birational contractions on Calabi--Yau threefolds of Type III.

\begin{example}\rm
Let $f_0\colon X_0\rightarrow \bar X_0$ be a birational 
contraction contracting a 
$\dP^1$-bundle $E_0$ to a curve $C_0$ of genus at least two. 
On a dense open subset $S$ of the marked moduli space of $X$, 
the contraction extends to a diagram 
\[\begin{array}{cccc}f\colon & \X & \rightarrow & \bar\X \\
& \downarrow & \swarrow  \\
& S.
\end{array}\]
The main new feature is that $f_s$ remains divisorial only for $s$ in a proper 
submanifold $S^{\rm div}\subset S$. At a general point $s\in S$, the 
contraction 
$f_s\colon X_s\rightarrow \bar X_s$ is of flopping type, i.e. it contracts 
a codimension two locus; I assume for simplicity that it contracts 
disjoint $\dP^1$s to nodes. 

By~\cite[Theorem 11.10]{kollar_mori}, there is a flop in the family 
\[\begin{array}{ccccc}\X && \dashrightarrow && \X^+ \\
& \searrow& & \swarrow\\
&& S
\end{array}\]
For $s\in S\setminus S^{\rm div}$ where $f_s$ is of flopping type, 
$X_s\dashrightarrow X_s^+$ is a classical flop. For $s\in S^{\rm div}$, 
the map $X_s\dashrightarrow X_s^+$ is the identity morphism. 

The crucial point is that $\X^+\rightarrow S$ is the pullback of the original 
family $\X\rightarrow S$ under a nontrivial map $\phi\colon S\rightarrow S$. 
The fixed locus of $\phi$ is exactly $S^{\rm div}$.  
There is an embedding \[\X\times_{\bar X}\X^+\hookrightarrow \X\times_S\X^+\]
and by definition, the latter space is isomorphic to $\X\times_\phi\X$. 
Thus the structure sheaf of $\X\times_{\bar X}\X^+$ gives an 
element $U\in\db{\X\times_\phi\X}$. 
\end{example}

\begin{theorem} The triple $(\phi, 0, U)$ gives an element of $G(\X,\M_X)$.  
\label{thm4}\end{theorem} 
{\vspace{.05in}\noindent {\sc Sketch Proof} \hspace{.05in}}
An easy argument shows that to prove that $(\phi, 0, U)$ is invertible, 
it is enough to show that the restriction is invertible on fibres. 
For $s\in S\setminus S^{\rm div}$, $U$ restricts 
to $X_s\times X_s^+$ as the structure sheaf of $X_s\times_{\bar X_s}X_s^+$. 
This is the well-known kernel of 
Bondal and Orlov~\cite[Theorem 3.6]{bon_orlov}, 
giving the isomorphism of derived categories under a classical flop. 

For points $s\in S^{\rm div}$, the restriction of $U$ 
to $X_s\times X_s$ is the kernel $\sheafO_{X_s\times_{\bar X_s}X_s}$.   
The method of Bridgeland~\cite{tom_flops} can be used to show that this 
sheaf is a universal sheaf for a certain moduli problem on $X_s$ with fine
moduli space $X_s$, and hence is an invertible kernel on $X_s\times X_s$. 
\ethrm 

\begin{remark} \rm 
The existence of the above kernel arising from Type III contractions was
conjectured by Horja in~\cite{horja}, although he used a different
expression~\cite[(4.35)]{horja}. 
In~\cite{horja2} he proves independently and 
by a different method that it is invertible. 
However, he does not discuss the relation to flops. 
\end{remark}

I illustrate the above discussion by a well-studied toric family. 

\begin{example} \label{typeiiiex} 
\rm
Start with a general degree eight hypersurface 
\[\bar X_0\subset\dP^4[1^2, 2^3]\] 
in weighted projective four-space. Standard toric technology or a simple blowup
give a crepant resolution $f_0\colon X_0\rightarrow\bar X_0$, which is a Type III 
contraction contracting a $\dP^1$-bundle $E_0$ in $X_0$ to a genus $3$ 
curve $C_0$ in $\bar X_0$. All toric deformations are resolutions of 
degree eight hypersurfaces in the weighted projective space; however, 
$h^{2,1}(X_0)=86$ whereas the octic has only 83 deformations, 
so some deformations of $X_0$ are missing. 

To find the other deformations, embed $\dP^4[1^2, 2^3]$ 
into $\dP^5$ using $\sheafO(2)$
as a hypersurface given by the equation $\{z_1z_2= z_3^2\}$. 
The variety $\bar X_0$ becomes a complete intersection
\setlength{\arraycolsep}{2pt}
\[\bar X_0\cong\left\{\begin{array}{lcl}z_1z_2& =& z_3^2 \\ f_4(z_i)&=&0\end{array}\right\}\subset\dP^5,\]
of a general quartic and a quadric of rank three in 
$\dP^5$. The curve of singularities is along $\{z_1=z_2=z_3=0\}$. 
Here more deformations of $X_0$ are visible:  
I can deform either the quadric or the quartic. 
Deformations of $f_4$ contribute 83 moduli, which are exactly the toric 
deformations; there are 3
further deformation directions deforming the quadric to a quadric of 
rank four\footnote{Deformations of $\bar X_0$ to non-singular complete intersections do not lift to deformations of the crepant resolution $X_0$.}. 
The generic such deformation 
\[\bar X_t=\left\{\begin{array}{lcl}z_1z_2&=&z_3^2-t^2z_4^2\\ f_4(z_i)&=&0\end{array}\right\}\subset\dP^5\]
has four nodes at $\{z_1=z_2=z_3=z_4=0\}$ and
a pair of small resolutions $X^\pm_t$, which are deformations of $X_0$:
\[ X_t^{\pm} = \left\{\begin{array}{lcl}z_1u_1 & = & (z_3\pm tz_4)u_2 \\(z_3\mp tz_4)u_1& = & z_2u_2 \\ f_4(z_i)& = & 0 \end{array}\right\}\subset\dP^1\times\dP^5.
\]
For $t\neq 0$, the birational map 
\[\psi_t\colon X_t^+\dashrightarrow X_t^-\] 
is a simple flop, flopping four copies
of $\dP^1$ over the four nodes of $\bar X_t$. 
By the main result of~\cite{tcex}, for general $t$ the
varieties $X_t^\pm$ are not isomorphic; what is more important from the point
of view of the present discussion is that the flopped family
explicitly appears as the original family under the base change 
$\phi: t\mapsto -t$. 

In the connected component $\X\rightarrow\M_X$ of the Calabi--Yau moduli
space of $X_0$, there is a sublocus $\M_X^{\rm div}\subset \M_X$
where the contraction $f_0$ deforms as a divisorial contraction; this is
also the fixed locus of the base change map $\phi$. 
On the other hand, as discussed above,
$\M_X^{\rm div}$ also equals $\M_X^{\rm toric}$, the toric deformation 
space of $X_0$. The element $(\phi, 0, U)$ of $G(\X, \M_X)$ guaranteed 
by Theorem~\ref{thm4} acts as a self-equivalence on the toric deformations
of $X_0$, but acts by Fourier--Mukai transforms between different
(generically non-isomorphic) deformations away from the toric locus. 
This is the same behaviour as that found at the end of 
Section~\ref{toricK3}.
\end{example} 

\subsection{Relation to the global Torelli problem} 

The global Torelli theorem for K3 surfaces is an important tool in the
proof of the cohomological form of Conjecture~\ref{mainconjecture} for K3
surfaces. For Calabi--Yau threefolds, 
the situation appears to be exactly the opposite. As 
an invertible kernel induces a Hodge isometry on $H^3$, the following holds:

\begin{proposition} 
Suppose that $\X\rightarrow \M_X$ is a component of the marked moduli space 
of a Calabi--Yau threefold. Assume that there is an element 
$(\phi, \alpha, U)\in G(\X,\M_X)$ such 
that for general $s\in \M_X$, the varieties $X_s$ and $X_{\phi(s)}$ 
are not isomorphic.
(In particular, $\phi$ is not the identity.) Then the global Torelli 
theorem fails for the family: the polarized 
(rational) Hodge structure on middle cohomology does not determine 
the general variety in the family up to isomorphism.  
\end{proposition}

The counterexample to global Torelli given in~\cite{tcex} indeed comes
from the nontrivial element of $G(\X,\M_X)$ of Theorem~\ref{thm4}. 
Some conjectural examples of Calabi--Yau threefold families with
a family of Fourier--Mukai transforms with nontrivial $\phi$ 
have been constructed by Andrei 
C\u al\-d\u a\-ra\-ru~\cite[Chapter 6]{andrei}. 
These examples are of a very different flavour; they are elliptic fibrations
without a section, the Fourier--Mukai transform coming from changing to
a different elliptic fibration with the same Jacobian.  
The framework presented here suggests that counterexamples to global Torelli 
arising from families of Fourier--Mukai transforms 
should be rather more common than hitherto suspected. 
However, one may speculate that the isomorphism of 
the transcendental part of Hodge structures (in this case simply $H^3$) may 
imply the isomorphism of derived categories; this would be at the end of 
the day in complete agreement with the K3 case. 
At least for cohomology with integral coefficients, 
there appears to be no known counterexample to this speculation.

\vspace{0.2in}

\noindent {\small Current address: 
Math. Inst., University of Warwick, Coventry CV4 7AL, UK} 

\noindent {\small Email: \tt balazs@maths.warwick.ac.uk} 

\end{document}